\documentclass[10pt,twoside]{article}
\usepackage{graphicx}
\usepackage{amsmath}
\usepackage{amssymb}
\usepackage[matrix,arrow,curve,cmtip]{xy}
\usepackage{Latex-document}

\newtheorem{theorem}[equation]{Theorem}
\newtheorem{conjecture}[equation]{Conjecture}
\newcommand{\xto}{\xrightarrow}
\renewcommand{\:}{\colon}
\newcommand{\Z}{\mathbb{Z}}
\newcommand{\Zp}{\mathbb{Z}_p}
\newcommand{\Fp}{\mathbb{F}_p}
\newcommand{\Qp}{\mathbb{Q}_p}
\newcommand{\s}{_{\boldsymbol{\cdot}}}
\newcommand{\cs}{^{\,\boldsymbol{\cdot}}}
\newcommand{\TR}{\operatorname{TR}}
\newcommand{\THH}{\operatorname{THH}}
\newcommand{\Spec}{\operatorname{Spec}}

\newcommand{\tate}{\hat{\mathbb{H}}}
\newcommand{\borel}{\mathbb{H}\hskip1pt_{\boldsymbol{\cdot}}}
\newcommand{\coborel}{\mathbb{H}^{\boldsymbol{\cdot}}}
\newcommand{\e}{\epsilon}
\newcommand{\f}{^{^{\wedge}}}

\title{\bf Algebraic {\boldmath $K$}-theory and Trace Invariants\vskip 6mm}

\author{Lars Hesselholt\vspace*{-0.5cm}\thanks{Massachusetts Institute
of Technology, Cambridge, MA 02139, USA. E-mail:
larsh@math.mit.edu}}

\date{\vspace{-8mm}}

\begin{document}

\maketitle

\begin{center}
\textit{$($Dedicated to Ib Madsen on his sixtieth birthday$)$}
\end{center}

\thispagestyle{first} \setcounter{page}{415}

\begin{abstract}

\vskip 3mm

The cyclotomic trace of B\"okstedt-Hsiang-Madsen, the subject of \linebreak B\"okstedt's lecture at the congress
in Kyoto, is a map of pro-abelian groups
$$K_*(A) \xrightarrow{\operatorname{tr}}
\operatorname{TR}_*^{\boldsymbol{\cdot}}(A;p)$$ from Quillen's algebraic $K$-theory to a topological refinement of
Connes' cyclic homology. Over the last decade, our understanding of the target and its relation to $K$-theory has
been significantly advanced. This and possible future development is the topic of my lecture.

The cyclotomic trace takes values in the subset fixed by an operator $F$ called the Frobenius. It is known that
the induced map
$$K_*(A,\mathbb{Z}/p^v) \xrightarrow{\operatorname{tr}}
\operatorname{TR}_*^{\boldsymbol{\cdot}}(A;p,\mathbb{Z}/p^v)^{F=1}$$ is an isomorphism, for instance, if $A$ is a
regular local $\mathbb{F}_p$-algebra, or if $A$ is a henselian discrete valuation ring of mixed characteristic
$(0,p)$ with a separably closed residue field. It is possible to evaluate $K$-theory by means of the cyclotomic
trace for a wider class of rings, but the precise connection becomes slightly more complicated to spell out.

The pro-abelian groups $\operatorname{TR}_*^{\boldsymbol{\cdot}}(A;p)$ are typically very large. But they come
equipped with a number of operators, and the combined algebraic structure is quite rigid. There is a universal
example of this structure --- the de~Rham-Witt complex
--- which was first considered by Bloch-Deligne-Illusie in
connection with Grothendieck's crystalline cohomology. In general, the canonical map
$$W_{\boldsymbol{\cdot}}\,\Omega_A^q \rightarrow
\operatorname{TR}_q^{\boldsymbol{\cdot}}(A;p)$$ is an isomorphism, if $q\leq 1$, and the higher groups, too, can
often be expressed in terms of the de~Rham-Witt groups. This is true, for example, if $A$ is a regular
$\mathbb{F}_p$-algebra, or if $A$ is a smooth algebra over the ring of integers in a local number field. The
calculation in the latter case verifies the Lichtenbaum-Quillen conjecture for local number fields, or more
generally, for henselian discrete valuation fields of geometric type.

\vskip 4.5mm


\end{abstract}

\vskip 12mm

\section{Algebraic {\boldmath $K$}-theory}\label{section 1}\setzero
\vskip-5mm \hspace{5mm}

The algebraic $K$-theory of Quillen~\cite{quillen}, inherently, is a
multiplicative theory. Trace invariants allow the study of this theory
by embedding it in an additive theory. It is possible, by this
approach, to evaluate the $K$-theory (with coefficients) of henselian
discrete valuation fields of mixed characteristic. We first recall the
expected value of the $K$-groups of a field $k$.

The groups $K_*(k)$ form a connected anti-commutative graded ring,
there is a canonical isomorphism $\ell\:k^*\xto{\sim}K_1(k)$, and
$\ell(x)\cdot\ell(1-x)=0$. One defines the Milnor $K$-groups
$K_*^M(k)$ to be the universal example of this algebraic
structure~\cite{milnor}. The canonical map
$K_q^M(k) \to K_q(k)$
is an isomorphism, if $q\leq 2$. Let us now fix the attention on the
$K$-groups with finite coefficients. (The rational $K$-groups,
while of great interest, are of a rather different
nature~\cite{dupont1,dupontsah}.) The groups $K_*(k,\Z/m)$ form an
anti-commutative graded $\Z/m$-algebra, at least if $v_2(m)\neq 1,2$
and $v_3(m)\neq 1$. And if $\mu_m\subset k$, there is a canonical
lifting
$$\xymatrix{
{} &
{K_2(k,\Z/m)} \ar[d]^{\beta_m} \cr
{\mu_m} \ar[r]^(.47){\ell} \ar@{-->}[ur]^{b} &
{K_1(k),} \cr
}$$
which to a primitive $m$th root of unity $\zeta$ associates the Bott
element $b_{\zeta}$. Hence, in this case, there is an additional map
of graded rings
$S_{\Z/m}(\mu_m) \to K_*(k,\Z/m).$
The Beilinson-Licthenbaum conjectures predict that the combined map
$$K_*^M(k)\otimes_{\Z}S_{\Z/m}(\mu_m) \to K_*(k,\Z/m)$$
be an isomorphism of graded rings~\cite{beilinson,licthenbaum1}. The
case $m=2^v$ follows from the celebrated proof of the Milnor
conjecture by Voevodsky~\cite{voevodsky}. We here consider the case of
a henselian discrete valuation field of mixed characteristic $(0,p)$
with $p$ odd and $m=p^v$~\cite{hm2,gh2}. The groups $K_*^M(k)/m$
typically are non-zero in only finitely many degrees. Hence, above
this range, the groups $K_*(k,\Z/m)$ are two-periodic. All rings
(resp.~graded rings, resp.~monoids) considered in this paper are
assumed commutative (resp.~anti-commutative, resp.~commutative) and
unital without further notice.

\section{The de~Rham-Witt complex}\label{section 2}\setzero
\vskip-5mm \hspace{5mm}

Let $V$ be a henselian discrete
valuation ring with quotient field $K$ of characteristic zero and
residue field $k$ of odd characteristic $p$. (At this writing, we
further require that $V$ be of geometric type, i.e. that $V$ be
the henselian local ring at the generic point of the special fiber
of a smooth scheme over a henselian discrete valuation ring
$V_0\subset V$ with \emph{perfect} residue field.) A first example of
a trace map is provided by the logarithmic derivative
$$K_*^M(K) \to \Omega_{(V,M)}^*$$
which to the symbol $\{a_1,\dots,a_q\}$ associates the form $d\log
a_1\dots d\log a_q$. The right hand side is the de~Rham complex
with log poles in the sense of Kato~\cite{kato}: A log ring $(A,M)$ is a
ring $A$ and a map of monoids $\alpha\:M\to (A,\,\cdot\,)$; a
log differential graded ring $(E^*,M)$ is a differential graded ring
$E^*$ together with maps of monoids $\alpha\:M\to (E^0,\cdot)$ and
$d\log\:M\to (E^1,+)$ such that $d\circ d\log=0$ and such that
$d\alpha(a)=\alpha(a)d\log a$ for all $a\in M$; the de~Rham complex
$\Omega_{(A,M)}^*$ is  the universal log differential graded ring with
underlying log ring $(A,M)$. We will always consider the ring $V$ with
the \emph{canonical} log structure
$$\alpha\:M=V\cap K^* \hookrightarrow V.$$
(In this case, there are natural short-exact sequences
$$0 \to \Omega_V^q \to \Omega_{(V,M)}^q \to \Omega_k^{q-1} \to 0.)$$
The logarithmic derivative, however, is far from injective. It turns
out that this can be rectified by incorborating the Witt vector
construction which we now recall.

The ring of Witt vectors associated with a ring $A$ is the set of
``vectors''
$$W(A)=\{(a_0,a_1,\dots)\mid a_i\in A\}$$
with a new ring structure, see~\cite{hm3}. The ring operations are
polynomial in the coordinates. The projection $W(A)\to A$, which to
$(a_0,a_1,\dots)$ associates $a_0$, is a natural ring homomorphism
with the unique natural multiplicative section
$$[\phantom{x}]\:A \to W(A), \hskip5mm [a]=(a,0,0,\dots).$$
If $F$ is a perfect field of characteristic $p>0$, then $W(F)$ is the
unique (up to unique isomorphism) complete discrete valuation ring
of mixed characteristic $(0,p)$ such that $W(F)/p\xto{\sim}F$. In
general, the ring $W(A)$ is equal to the inverse limit of the rings
$W_n(A)$ of Witt vectors of length $n$. But rather than forming the
limit, we shall consider the limit system of rings $W\s(A)$ as a
pro-ring. There is a natural map of pro-rings $F\: W\s(A) \to
W_{\boldsymbol{\cdot}-1}(A)$, called the Frobenius, and a natural map
of $W\s(A)$-modules $V\: F_*\,W_{\boldsymbol{\cdot}-1}(A) \to W\s(A)$,
called the Verschiebung. The former is given, as the ring structure,
by certain polynomials in the coordinates; the latter is given by
$V(a_0,\dots,a_{n-2})=(0,a_0,\dots,a_{n-2})$, and $FV=p$. Finally, we
note that if $(A,M)$ is a log ring, then the composite
$$M \xto{\alpha} A \xto{[\phantom{a}]} W\s(A)$$
makes $(W\s(A),M)$ a pro-log ring.

There is a natural way to combine differential forms and Witt
vectors; the result is called the de~Rham-Witt complex. It was
considered first for $\Fp$-algebras by
Bloch-Deligne-Illusie~\cite{bloch,illusie} in  connection with the
crystalline cohomology of Berthelot-Grothendieck~\cite{berthelot}. A
generalization to log-$\Fp$-algebras was constructed by
Hyodo-Kato~\cite{hyodokato}. The following extension to
log-$\Z_{(p)}$-algebras was obtained in collaboration with Ib
Madsen~\cite{hm3,hm2}: Let $(A,M)$ be a log ring such that $A$ is a
$\Z_{(p)}$-algebra with $p$ odd. A \emph{Witt complex} over $(A,M)$
is:
\vskip1mm

(i) a pro-log differential graded ring $(E\s^*,M_E)$ and a map of
pro-log rings
$$\lambda\:(W\s(A),M)\to (E\s^0,M_E);$$

(ii) a map of pro-log graded rings
$$F\:E\s^*\to E_{\boldsymbol{\cdot}-1}^*,$$
such that $\lambda F=F\lambda$ and such that
$$\begin{aligned}
Fd\log_n a & = d\log_{n-1}a, \hskip15.2mm
\text{for all $a\in M$,} \cr
Fd\,\lambda[a]_n & =\lambda[a]_{n-1}^{p-1}d\,\lambda[a]_{n-1},
\hskip5mm \text{for all $a\in A$;}\cr
\end{aligned}$$

(iii) a map of pro-graded modules over the pro-graded ring $E^*\s$,
$$V\:F_*\,E_{\boldsymbol{\cdot}-1}^*\to E\s^*,$$
such that $\lambda V=V\lambda$, $FV=p$ and $FdV=d$.

A map of Witt complexes over $(A,M)$ is a map of pro-log differential
graded rings which commutes with the maps $\lambda$, $F$ and $V$.
Standard category theory shows that there exists a universal Witt
complex over $(A,M)$. This, by defintion, is the de~Rham-Witt complex
$W\s\,\Omega_{(A,M)}^*$. (The canonical maps $W\s(A) \to
W\s\,\Omega_{(A,M)}^0$ and $\Omega_{(A,M)}^* \to W_1\,\Omega_{(A,M)}^*$
are isomorphisms, so the construction really does combine differential
forms and Witt vectors.) We  lift the logarithmic derivative to
a map
$$K_q^M(K) \to W_n\,\Omega_{(V,M)}^q$$
which to the symbol $\{a_1,\dots,a_q\}$ associates $d\log_na_1\dots
d\log_na_q$. This trace map better captures the Milnor $K$-groups.
Indeed, the following result was obtained in collaboration with
Thomas Geisser~\cite{gh2}:

\begin{theorem}\label{kdrw} Suppose that $\mu_{p^v}\subset K$ and that
$k$ is separably closed. Then the trace map induces an isomorphism of
pro-abelian groups
$$K_q^M(K)/p^v \xto{\sim}
\big(W\s\,\Omega_{(V,M)}^q/p^v\big)^{F=1}.$$
\end{theorem}

To prove this, we first show that $W_n\,\Omega_{(V,M)}^q/p$ has a
(non-canonical) $k$-vector space structure and find an explicit
basis. The dimension is
$$\operatorname{dim}_k\big(W_n\,\Omega_{(V,M)}^q/p\big) =
\binom{r+1}{q}\,e\,\sum_{s=0}^{n-1}p^{rs},$$
where $|k:k^p|=p^r$ and $e$ the ramification index of $K$. It is not
difficult to see that this is an upper bound for the dimension. The
proof that it is also a lower bound is more involved and uses a
formula for the de~Rham-Witt complex of a polynomial extension by
Madsen and the author~\cite{hm3}. We then evaluate the kernel of $1-F$
and compare with the calculation of $K_q^M(K)/p$ by
Kato~\cite{kato1,blochkato}. The assumption that the residue field $k$
be separably closed is not essential. In the general case, one instead
has a short-exact sequence
$$0 \to \big(W\s\,\Omega_{(V,M)}^{q-1}\otimes\mu_{p^v}\big)_{F=1}
\to K_q^M(K)/p^v \to \big(W\s\,\Omega_{(V,M)}^q/p^v\big)^{F=1} \to
0,$$
where the superscript (resp.~subscript) ``$F=1$'' indicates Frobenius
invariants (resp.~coinvariants).

We discuss a global version of theorem~\ref{kdrw}. Let $V_0$ be a
henselian discrete valuation ring with quotient field $K_0$ of
characteristic zero and \emph{perfect} residue field $k_0$ of odd
characteristic $p$. Let $\mathfrak{X}$ be a smooth $V_0$-scheme, and
let $i$ (resp. $j$) denote the inclusion of the special
(resp. generic) fiber as in the cartesian diagram
$$\xymatrix{
{X\;} \ar@{^{(}->}[r]^{j} \ar[d] &
{\mathfrak{X}} \ar[d]^{f} &
{\;Y} \ar@{_{(}->}[l]_{i} \ar[d] \cr
{\Spec K_0\;} \ar@{^{(}->}[r] &
{\Spec V_0} &
{\;\Spec k_0.} \ar@{_{(}->}[l] \cr
}$$
Suppose that $\mu_{p^v}\subset K_0$. Then the proof of
theorem~\ref{kdrw} shows that there is a short-exact sequence of
sheaves of pro-abelian groups on $Y$ for the \'etale topology
$$0 \to i^*R^qj_*\Z/p^v(q) \to
i^*(W\s\,\Omega_{(\mathfrak{X},M)}^q/p^v) \xto{1-F}
i^*(W\s\,\Omega_{(\mathfrak{X},M)}^q/p^v) \to 0.$$
The left hand term is the sheaf of $p$-adic vanishing cycles.

\section{The cyclotomic trace}\label{section 3}\setzero
\vskip-5mm \hspace{5mm}

We now turn to Quillen $K$-theory. The
analog of the logarithmic derivative is the topological Dennis trace
with values in topological Hochschild homology,
$$K_*(\mathcal{C}) \to \THH_*(\mathcal{C}),$$
defined by B\"okstedt~\cite{bokstedt}. It is a refinement of earlier
trace maps by Dennis~\cite{dennis} and Waldhausen~\cite{waldhausen1}.
We will use a variant of the construction due to
Dundas-McCarthy~\cite{dundasmccarthy,mccarthy} that can be
applied to a category with cofibrations and weak equivalences in the
sense of Waldhausen~\cite{waldhausen}. The
category $\mathcal{C}$ we consider is the category of bounded chain
complexes of finitely generated projective $V$-modules. The
cofibrations are the degree-wise monomorphisms, and the weak
equivalences are the chain maps $C\to C'$ such that
$K\otimes_VC\to K\otimes_VC'$ is a quasi-isomorphism. The
$K$-theory of this category is canonically isomorphic to Quillen's
$K$-theory of the field $K$. We showed in~\cite{hm2} that the groups
$$\THH_*(V|K)=\THH_*(\mathcal{C})$$
form a log differential graded ring with underlying log ring
$(V,M)$, where the structure map $d\log$ is given by the composite
$$M=V\cap K^* \xto{\ell} K_1(K) \to \THH_1(V|K).$$
The canonical map from the de~Rham complex
$$\Omega_{(V,M)}^q \to \THH_q(V|K)$$
is compatible with the trace maps and is an isomorphism, if $q\leq 2$.
The topological Dennis trace, again, is far from injective. This can be
rectified by a construction which, in retrospect, can be seen as
incorporating Witt vectors. The result is the cyclotomic trace of
B\"okstedt-Hsiang-Madsen~\cite{bhm} which we now recall. The reader is
referred to~\cite{hm2,gh,dundasmccarthy} for details.

The topological Dennis trace, we recall, is defined as the map of
homotopy groups induced from a continuous map of spaces
$$K(\mathcal{C}) \to \THH(\mathcal{C}).$$
As a consequence of Connes' theory of cyclic sets, the right hand
space is equipped with a continuous action by the circle group
$\mathbb{T}$. Moreover, the image of the trace map is point-wise fixed
by the $\mathbb{T}$-action. Let
$$\TR^n(\mathcal{C};p)=\THH(\mathcal{C})^{C_{p^{n-1}}}$$
be the subset fixed by the subgroup $C_{p^{n-1}}\subset\mathbb{T}$ of
the indicated order. It turns out that, as $n$ and $q$ varies, the
homotopy groups
$$\TR_q^n(V|K;p)=\pi_q(\TR^n(\mathcal{C};p))$$
form a Witt complex over $(V,M)$; see~\cite{hm,h,hm2}. The map $F$ is
induced from the obvious inclusion map, and the map $V$ is the
accompanying transfer map. The structure maps in the limit system and
the map $\lambda$, however, are more difficult to define. The former
was defined in~\cite{bhm} and the latter in~\cite{hm}. The topological
Dennis trace induces a map of pro-abelian groups
$$K_q(K) \to \TR_q\cs(V|K;p).$$
This is the cyclotomic trace. It takes values in the subset fixed
by the Frobenius operator. The canonical map
$$W_n\Omega_{(V,M)}^q \to \TR_q^n(V|K;p)$$
is compatible with the trace maps from Milnor $K$-theory and Quillen
$K$-theory, respectively, and is an isomorphism, if $q\leq 2$. The
following is a combination of results obtained in collaboration with
Thomas Geisser~\cite{gh} and Ib Madsen~\cite{hm,hm2}.

\begin{theorem}\label{trace}Suppose that $k$ is separably closed. Then
the cyclotomic trace induces an isomorphism of pro-abelian groups
$$K_q(K,\Z/p^v) \xto{\sim} \TR_q\cs(V|K;p,\Z/p^v)^{F=1}.$$
\end{theorem}

We briefly outline the steps in the proof: We proved in~\cite{gh} that
the sequence
$$0 \to K_q(k,\Z/p^v) \to \TR_q\cs(k;p,\Z/p^v) \xto{1-F}
\TR_q\cs(k;p,\Z/p^v) \to 0$$
is exact. This uses~\cite{blochkato, geisserlevine,h}. Given this, the
theorem by McCarthy~\cite{mccarthy1} that for nilpotent extensions,
relative $K$-theory and relative topological cyclic homology agree,
and the continuity results of Suslin~\cite{suslin1} for $K$-theory and
Madsen and the author~\cite{hm} for $\TR$ show that also the sequence
$$0 \to K_q(V,\Z/p^v) \to \TR_q\cs(V;p,\Z/p^v) \xto{1-F}
\TR_q\cs(V;p,\Z/p^v) \to 0$$
is exact. Theorem~\ref{trace} follows by comparing the localization
sequence of Quillen~\cite{quillen}
$$\cdots\to K_q(k,\Z/p^v) \xto{i^!} K_q(V,\Z/p^v) \xto{j_*}
K_q(K,\Z/p^v) \to \cdots$$
to the corresponding sequence by Madsen and the author~\cite{hm2}
$$\cdots\to \TR_q^n(k;p,\Z/p^v) \xto{i^!} \TR_q^n(V;p,\Z/p^v)
\xto{j_*} \TR_q(V|K;p,\Z/p^v) \to \cdots.$$
Again, the assumption in the statement of theorem~\ref{trace} that the
residue field $k$ be separably closed is not essential. The general
statement will be given below. It is also not necessary for
theorem~\ref{trace} to assume that $V$ be of geometric type.

\section{The Tate spectral sequence}\label{section 4}\setzero
\vskip-5mm \hspace{5mm}

If $G$ is a finite group and $X$ a
$G$-space, it is usually not possible to evaluate the groups
$\pi_*(X^G)$ from knowledge of the $G$-modules $\pi_*(X)$. At first
glance, this is the problem that one faces in evaluating the groups
$$\TR_q^n(\mathcal{C};p) =
\pi_q\big(\THH(\mathcal{C})^{C_{p^{n-1}}}\big).$$
However, the mapping fiber of the structure map
$\TR^n(\mathcal{C};p) \to \TR^{n-1}(\mathcal{C};p)$, it turns out,
is given by the Borel construction
$\borel(C_{p^{n-1}},\THH(\mathcal{C}))$ whose homotopy groups are the
abutment of a (first quadrant) spectral sequence
$$E_{s,t}^2=H_s(C_{p^{n-1}},\THH_t(\mathcal{C}))
\Rightarrow \pi_{s+t}\borel(C_{p^{n-1}},\THH(\mathcal{C})).$$
This suggests that the groups $\TR_q^n(\mathcal{C};p)$ can be
evaluated inductively starting from the case $n=1$. However, it is
generally difficult to carry out the induction step. In addition, the
absence of a multiplicative structure makes the spectral sequence
above difficult to solve. The main vehicle to overcome these
problems, first employed by B\"okstedt-Madsen
in~\cite{bokstedtmadsen}, is the following diagram of fiber sequences
$$\xymatrix{
{\borel(C_{p^{n-1}},\THH(\mathcal{C}))} \ar[r] \ar@{=}[d] &
{\TR^n(\mathcal{C};p)} \ar[r] \ar[d]^{\Gamma} &
{\TR^{n-1}(\mathcal{C};p)} \ar[d]^{\hat{\Gamma}} \cr
{\borel(C_{p^{n-1}},\THH(\mathcal{C}))} \ar[r] &
{\coborel(C_{p^{n-1}},\THH(\mathcal{C}))} \ar[r] &
{\tate(C_{p^{n-1}},\THH(\mathcal{C}))} \cr
}$$
together with a multiplicative (upper half-plane) spectral sequence
$$E_{s,t}^2 = \hat{H}^{-s}(C_{p^{n-1}},\THH_t(\mathcal{C}))
\Rightarrow \pi_{s+t}\tate(C_{p^{n-1}},\THH(\mathcal{C}))$$
starting from the Tate cohomology of the (trivial)
$C_{p^{n-1}}$-module $\THH_t(\mathcal{C})$. The lower fiber sequence
is the Tate sequence; see Greenlees and May~\cite{greenleesmay}
or~\cite{hm2}. In favorable cases, the maps $\Gamma$ and
$\hat{\Gamma}$ induce isomorphisms of homotopy groups in non-negative
degrees. Indeed, this is true in the case at hand (if $k$ is
perfect). The differential structure of the spectral sequence
$$E_{s,t}^2 = \hat{H}^{-s}(C_{p^{n-1}},\THH_t(V|K,\Z/p))
\Rightarrow \pi_{s+t}(\tate(C_{p^{n-1}},\THH(V|K)),\Z/p)$$
was determined in collaboration with Ib Madsen~\cite{hm2} in the
case where the residue field $k$ is perfect. This is the main
calculational result of the work reported here. The following
result, for perfect $k$, is a rather immediate consequence. The
extension to non-perfect $k$ is given in~\cite{hm3}.

\begin{theorem}\label{main}Suppose that $\mu_{p^v}\subset K$. Then the
canonical map is an isomorphism of pro-abelian groups
$$W\s\,\Omega_{(V,M)}^*\otimes_{\Z}S_{\Z/p^v}(\mu_{p^v}) \xto{\sim}
\TR_*\cs(V|K;p,\Z/p^v).$$
\end{theorem}

We can now state the general version of theorem~\ref{trace} which does
not require that the residue field $k$ be separably closed. The second
tensor factor on the left hand side in the statement of
theorem~\ref{main} is the symmetric algebra on the $\Z/p^v$-module
$\mu_{p^v}$, which is free of rank one. Spelling out the statement for
the group in degree $q$, we get an isomorphism of pro-abelian groups
$$\bigoplus_{s\geq
0}W\s\,\Omega_{(V,M)}^{q-2s}\otimes\mu_{p^v}^{\otimes s}
\xto{\sim}\TR_q\cs(V|K;p,\Z/p^v).$$
In the case of a separably closed residue field, theorem~\ref{trace}
idenfies the Frobenius fixed set of the common pro-abelian group
with $K_q(K,\Z/p^v)$. In the general case, one has instead a
short-exact sequence
$$0 \to \bigoplus_{s\geq 1}
\big(W\s\,\Omega_{(V,M)}^{q+1-2s} \otimes
\mu_{p^v}^{\otimes s}\big)_{F=1}\!\! \to
K_q(K,\Z/p^v) \to
\bigoplus_{s\geq 0}
\big(W\s\,\Omega_{(V,M)}^{q-2s} \otimes
\mu_{p^v}^{\otimes s}\big)^{F=1}\!\! \to 0,$$
valid for all integers $q$. (There is a similar sequence for the
topological cyclic homology group
$\operatorname{TC}_q\cs(V|K;p,\Z/p^v)$~\cite{hm2} which includes the
summand ``$s=0$'' on the left.) Comparing with the general version of
theorem~\ref{kdrw}, we obtain the following result promised
earlier~\cite{hm2,gh2}.

\begin{theorem}\label{bl}Suppose that $\mu_{p^v}\subset K$. Then the
canonical map
$$K_*^M(K)\otimes_{\Z}S_{\Z/p^v}(\mu_{p^v}) \xto{\sim}
K_*(K,\Z/p^v)$$
is an isomorphism.
\end{theorem}

\section{Galois descent}\label{section 5}\setzero
\vskip-5mm \hspace{5mm}

We now assume that the residue field $k$ be perfect. In homotopy
theoretic terms, theorem~\ref{main} states that the pro-spectrum
$\TR\cs(V|K;p)$ is equivalent to the $(-1)$-connected cover of its
localization with respect to complex periodic $K$-theory,
see~\cite{bousfield}. This suggests the possibility
of completely understanding the homotopy type of this pro-spectrum. We
expect that this, in turn, is closely related to the following question.
Let $\bar{K}$ be an algebraic closure of $K$ with Galois group $G_K$,
and let $\bar{V}$ be the integral closure of $V$ in $\bar{K}$. (The
ring $\bar{V}$ is a valuation ring with value group the additive group
of rational numbers.)

\begin{conjecture}If $k$ is perfect then for all $q>0$, the canonical map
$$\TR_q\cs(V|K;p,\Qp/\Zp) \to
\TR_q\cs(\bar{V}|\bar{K};p,\Qp/\Zp)^{G_K}$$
be an isomorphism of pro-abelian groups and that the higher continuous
cohomology groups
$H_{\operatorname{cont}}^i(G_K,\TR_q^n(\bar{V}|\bar{K};p,\Qp/\Zp))$
vanish.
\end{conjecture}

It follows from Tate~\cite{tate} that the groups
$H_{\operatorname{cont}}^i(G_K,\TR_q^n(\bar{V}|\bar{K};p,\Qp))$
vanish for $i\geq 0$ and $q>0$. One may hope that these methods will
help shed some light on the structure of the groups
$H_{\operatorname{cont}}^i(G_K,\TR_q^n(\bar{V}|\bar{K};p,\Qp/\Zp))$.
We now describe the structure of these $G_K$-modules; proofs will
appear elsewhere.

The group $\TR_q^n(\bar{V}|\bar{K};p,\Qp/\Zp)$ is divisible, if $q>0$,
and uniquely divisible, if $q>0$ and even. The Tate module
$T_p\TR_1^n(\bar{V}|\bar{K};p)$ is a free module of rank one over
$\TR_0^n(\bar{V}|\bar{K};p,\Zp)$, and the canonical map an
isomorphism:
$$S_{\TR_0^n(\bar{V}|\bar{K};p,\Zp)}(T_p\TR_1^n(\bar{V}|\bar{K};p))
\xto{\sim}\TR_*^n(\bar{V}|\bar{K};p,\Zp)$$
(note that  $\TR_q^n(\bar{V}|\bar{K};p,\Qp/\Zp) \xto{\sim}
\TR_q^n(\bar{V}|\bar{K};p,\Zp)\otimes\Qp/\Zp$). We note the formal
analogy with the results on $K_*(\bar{K})$ by
Suslin~\cite{suslin,suslin1}.

The structure of the ring
$\TR_0^n(\bar{V}|\bar{K};p,\Zp)=W_n(\bar{V})\f$ is well-understood
(unlike that of $W_n(V)$): Following Fontaine~\cite{fontaine}, we
let $R_{\bar{V}}$ be the inverse limit of the diagram $\bar{V}/p
\leftarrow \bar{V}/p \leftarrow \cdots$ with the Frobenius as
structure map. This is a perfect $\Fp$-algebra and an integrally
closed domain whose quotient field is algebraically closed. There is a
surjective ring homomorphism
$\theta_n\: W(R_{\bar{V}}) \twoheadrightarrow W_n(\bar{V})\f$
whose kernel is a principal ideal. If $\e=\{\e^{(v)}\}_{v\geq 1}$ is a
compatible sequence of primitive $p^{v-1}$st roots of unity considered
as an element of $R_{\bar{V}}$, and if $\e_n$ is the unique $p^n$th
root of $\e$, then $([\e]-1)/([\e_n]-1)$ is a generator. Moreover, as
$n$ varies, the maps $\theta_n$ constitute a map of pro-rings
compatible with the Frobenius maps.

The Bott element $b_{\e,n}\in T_p\TR_1^n(\bar{V}|\bar{K};p)$
determined by the sequence $\e$ is not a generator (so the statement
of theorem~\ref{main} is not valid for $\bar{K}$). Instead there is a
generator $\alpha_{\e,n}$ such that
$b_{\e,n}=([\e_n]-1)\cdot\alpha_{\e,n}$. The structure map of the
pro-abelian group $T_p\TR_1\cs(\bar{V}|\bar{K};p)$ (resp. the
Frobenius) takes $\alpha_{\e,n}$ to
$([\e_{n-1}]-1)/([\e_n]-1)\cdot\alpha_{\e,n-1}$ (resp. to
$\alpha_{\e,n-1}$), and the action of the Galois group is given by
$$\alpha_{\e,n}^{\sigma} =
\chi(\sigma)\frac{[\e_n]-1}{[\e_n^{\sigma}]-1}\cdot\alpha_{\e,n},$$
where $\chi\:G_K \to \operatorname{Aut}(\mu_{p^{\infty}})=\Zp^*$ is
the cyclotomic character.

\bigskip

\noindent{\bf Acknowledgments} The research reported here was
supported in part by grants from the National Science Foundation
and by an Alfred P.~Sloan Fellowship.

\providecommand{\bysame}{\leavevmode\hbox to3em{\hrulefill}\thinspace}

\label{lastpage}

\end{document}